\documentclass[12pt,a4paper]{article}
\usepackage[a4paper,margin=15mm]{geometry}
\usepackage{amsmath,amssymb,amsthm}
\usepackage[T1]{fontenc}
\usepackage{microtype}
\usepackage[hidelinks]{hyperref}
\hypersetup{pdftitle={An asymptotic form of Yuster's independent transversal factor conjecture},pdfauthor={Xin Zhang}}

\newtheorem{theorem}{Theorem}
\newtheorem{lemma}{Lemma}
\newtheorem{claim}{Claim}
\newtheorem{conjecture}{Conjecture}
\newcommand{\Prob}{\mathbb P}
\newcommand{\E}{\mathbb E}

\title{An asymptotic bound for Yuster's independent
transversal factor conjecture}
\author{Xin Zhang\\
School of Mathematics and Statistics, Xidian University\\
Xi'an 710071, China\\
\texttt{xzhang@xidian.edu.cn}}
\date{}

\begin{document}
\maketitle

\begin{abstract}
An \emph{$[r,n,1]$-graph} has $r$ independent parts of size $n$, with a matching between each pair of parts. Yuster conjectured in 2021 that every such graph has a partition into independent transversals whenever $n\ge 2\lceil r/2\rceil$. We prove that there is an absolute constant $A>0$ such that, for all sufficiently large $r$, this partition exists whenever $n\ge\lceil r+A r^{20/21}\rceil$. Thus $f(r)\le r+O(r^{20/21})$, where $f(r)$ is the least $n_0$ such that every $[r,n,1]$-graph has such a partition whenever $n\ge n_0$. In particular, $f(r)/r\to1$, giving an asymptotic form of the conjecture.
\end{abstract}

\section{Introduction}

Let $G$ have independent parts $V_1,\ldots,V_r$, each of size $n$.
The graph $G$ is an \emph{$[r,n,1]$-graph} if the edges between any two
distinct parts form a matching. These matchings may be partial.
An \emph{independent transversal} of $G$ is an independent set containing
exactly one vertex from each part. A \emph{factor of independent
transversals} is a partition of $V(G)$ into $n$ independent
transversals. Let $f(r)$ denote the least integer $n_0$ such that
every $[r,n,1]$-graph admits a factor of independent transversals
whenever $n\ge n_0$.

Yuster~\cite{Yuster97} studied factors of independent transversals in
sparse partite graphs. MacKeigan~\cite{MacKeigan} related these
factors to pairs of proper colorings such that distinct vertices
receive distinct pairs of colors. Hall's theorem gives
$f(r)\le 2r-2$; see~\cite{Yuster21}. Using random extensions, Yuster \cite{Yuster21}
proved that $f(r)\le 1.78r$ for all sufficiently large $r$.

For a lower bound, take $r-1$ disjoint copies of $K_r$, with each
copy meeting every part. An independent set contains at most one
vertex from each copy, so this $[r,r-1,1]$-graph has no factor of
independent transversals. Hence
\begin{equation}\label{eq:lower-bound}
 f(r)\ge r.
\end{equation}
For odd $r$, a
construction in~\cite{Yuster21} gives the stronger bound
$f(r)\ge r+1$.

Yuster proposed that these lower bounds give the exact value of $f(r)$.
\begin{conjecture}[Yuster~\cite{Yuster21}]\label{conj:yuster}
For every integer $r\ge2$,
\[
 f(r)=2\lceil r/2\rceil.
\]
\end{conjecture}

The \emph{local degree} of a multipartite graph is the maximum number of
neighbors that a vertex has in any one other part. Haxell~\cite{Haxell}
proved that a graph has an independent transversal if every part has
size at least $2\Delta$, where $\Delta$ is its maximum degree. For
each fixed $\varepsilon>0$, Loh and Sudakov~\cite{LohSudakov} lowered
the required part size to $(1+\varepsilon)\Delta$ when the local
degree is $o(\Delta)$. Glock and Sudakov~\cite{GlockSudakov}, and
independently Kang and Kelly~\cite{KangKelly}, showed that a bound
on the average vertex degree in each part can replace the
maximum-degree condition.

The \emph{multipartite complement} $\overline{G}$ of $G$ is obtained by
joining two vertices in different parts precisely when they are not
adjacent in $G$. A factor of independent transversals in $G$ is then a
perfect $K_r$-packing in $\overline{G}$.

Keevash and Mycroft~\cite{KeevashMycroft} proved a multipartite version
of the Hajnal--Szemer\'edi theorem. For fixed integers $R\ge k\ge 2$,
there is an $n_0=n_0(R,k)$ such that, whenever $n\ge n_0$ and
$k\mid Rn$, every balanced $R$-partite graph with parts of size $n$
has a perfect $K_k$-packing if each vertex has at least
$(k-1)n/k$ neighbors in every other part, unless
$Rn/k$ is odd, $k\mid n$, and the graph is isomorphic to
$\Gamma_{n,R,k}$. In this theorem, $R$ and $k$ are fixed while $n$ grows, and the
threshold $n_0$ may depend on both parameters. Thus, the theorem does
not provide a bound that is uniform when the number of parts and the
clique size grow with the part size.

The same problem can be stated in terms of \emph{correspondence packing}~\cite{PackingLists,BoundedPacking}. An \emph{$n$-fold correspondence cover} of $K_r$ has one independent set of $n$ color choices for each vertex of $K_r$. Between each pair of sets, a matching specifies the pairs of choices that cannot be used together. A \emph{correspondence coloring} selects one choice from each set and avoids all forbidden pairs. A \emph{full correspondence packing} consists of $n$ such colorings that use every choice exactly once. Thus these covers are precisely the $[r,n,1]$-graphs, and their full packings are factors of independent transversals.

The \emph{correspondence packing number} of $K_r$, also called its \emph{DP-packing number}, is denoted by $\chi_c^\star(K_r)$. It is the least $n$ for which every $n$-fold cover has a full packing. This property also holds for every larger cover. Indeed, restricting each choice set to $n$ elements gives a coloring. Remove such colorings until each set has size $n$, and then pack the remaining choices. Hence
\begin{equation}\label{eq:packing-number}
 f(r)=\chi_c^\star(K_r).
\end{equation}
Cambie, Cames van Batenburg, Davies and Kang~\cite{PackingLists} discuss Yuster's conjecture and his upper bound in this notation.

For packings, Chakraborti and Tran~\cite{ChakrabortiTran} proved that,
for each fixed $0<\beta<1$, a sufficiently small local degree
guarantees at least $\lfloor(1-\beta)n\rfloor$ disjoint independent
transversals when the parts have size $n$ and the maximum degree is at
most $(1-\beta)n$. Now fix $\varepsilon>0$ and set
$\beta=\varepsilon/(1+\varepsilon)$. If an $[r,n,1]$-graph satisfies
$n\ge(1+\varepsilon)r$, then its maximum degree is at most
\[
r-1\le \frac{n}{1+\varepsilon}=(1-\beta)n,
\]
and its local degree is at most $1$. Thus, for sufficiently large
$n$, their theorem gives at least $\lfloor(1-\beta)n\rfloor$
disjoint independent transversals, leaving at most
\[
\left\lceil\frac{\varepsilon n}{1+\varepsilon}\right\rceil
\]
vertices uncovered in each part.

The main result of this paper shows that every $[r,n,1]$-graph has a factor of
independent transversals when $n$ is at least $r+A r^{20/21}$,
for an absolute constant $A$ and all sufficiently large $r$.

\begin{theorem}\label{thm:main}
There are an absolute constant $A>0$ and an integer $r_0$ such that
every $[r,n,1]$-graph with $r\ge r_0$ and
\[
 n\ge \left\lceil r+A r^{20/21}\right\rceil
\]
has a factor of independent transversals.
\end{theorem}

Together with the lower bound $f(r)\ge r$ in~\eqref{eq:lower-bound}, Theorem~\ref{thm:main}
gives
\[
f(r)=r+O(r^{20/21}),
\qquad\text{and hence}\qquad
\frac{f(r)}{r}\longrightarrow 1.
\]
By~\eqref{eq:packing-number}, this also gives
\[
 \chi_c^\star(K_r)\le r+O(r^{20/21}).
\]

To prove the theorem, we build the factor one part at a time. At each
step, a bipartite graph records which vertices in the next part can
extend the current independent sets. We use edge weights
to choose a matching and then complete it to a perfect matching by
changing only a few edges. Throughout the process, we preserve enough
possible extensions to continue through all the parts.

Section~\ref{sec:setup} describes the construction. The matching
tools it uses are developed in Sections~\ref{sec:matching}
and~\ref{sec:random}. Section~\ref{sec:proof} shows that the
construction can be completed through all parts.

\section{Setup}\label{sec:setup}

Let $G$ be an $[r,n,1]$-graph with parts $V_1,\ldots,V_r$.
Fix $0<\rho<1$ and suppose that $2\le r\le \rho n$.
Assign distinct labels $1,\ldots,n$ to the vertices of $V_1$. At
stage $1$, label $i$ indexes the singleton containing the vertex of
$V_1$ with label $i$. In general, at stage $s\ge 2$, the labels index a partition of
$V_1\cup\cdots\cup V_s$ into independent sets, each containing one
vertex from every part $V_1,\ldots,V_s$. Call these sets
\emph{partial independent transversals}.
The \emph{compatibility graph} $B_s$ is the bipartite graph with
labels on the left and $V_{s+1}$ on the right, where $iv$ is an edge
if adding $v$ to the partial independent transversal indexed by $i$
preserves independence. More generally, a label and a vertex in any
part not yet assigned are called \emph{compatible} if adding that
vertex to the set indexed by the label preserves independence.
Each vertex on either side of $B_s$ has at most $s$ non-neighbors. For a left label this follows because each vertex of its transversal forbids at most one vertex in the next part. For a right vertex it follows because the vertex has at most one neighbor in each old part. A perfect matching in $B_s$ assigns every vertex of $V_{s+1}$ to a different label and extends all partial transversals.

To continue the construction, we impose a condition at each stage.
Choose parameters $1<a<1/\rho$ and $\xi>0$, and put
\begin{equation}\label{eq:stage-parameters}
 q=\lfloor \xi n\rfloor,\qquad
 p_s=1-\frac{as}{n},\qquad
 T_s=\frac{a}{p_s}.
\end{equation}
Before assigning part $V_{s+1}$, require that every set $Q$ of $q$
labels and every $W\subseteq V_{s+1}$ satisfy
\begin{equation}\label{eq:check}
 \frac{1}{qn}\sum_{i\in Q}\deg_W(i)\ge p_s\theta-\xi,
 \qquad \theta=\frac{|W|}{n}.
\end{equation}
Here $\deg_W(i)$ is the number of vertices in $W$ compatible with
label $i$. If there is a set $Q$ of $q$ labels and a set
$W\subseteq V_{s+1}$ for which~\eqref{eq:check} fails, stop the
process.

We now fix the parameters and record their bounds. Set
\begin{equation}\label{eq:parameter-choice}
 \alpha=1-\rho,\qquad a=1+\frac{\alpha}{2},\qquad
 \xi=\frac{\alpha^3}{100},\qquad
 p_*=1-a\rho,\qquad T_* = \frac{a}{p_*}.
\end{equation}

\begin{lemma}\label{lem:parameters}
For every $1\le s<r$, these parameters satisfy
\begin{equation}\label{eq:stage-bounds}
 1<a<1/\rho,\qquad p_s\ge p_*\ge\alpha/2,\qquad
 T_s\le T_*\le3/\alpha.
\end{equation}
They also satisfy
\begin{equation}\label{eq:slack-bounds}
 T_*\xi<(a-1)\alpha,\qquad \xi<p_*\alpha<\alpha.
\end{equation}
If $1/n<\xi/4$, then
\begin{equation}\label{eq:q-bounds}
 1\le \xi n/2\le q\le n.
\end{equation}
\end{lemma}

\begin{proof}
Since $0<\alpha=1-\rho<1$, we have $1<a<3/2$ and
\[
 1-a\rho
 =1-\left(1+\frac{\alpha}{2}\right)(1-\alpha)
 =\frac{\alpha(1+\alpha)}{2}\ge\frac{\alpha}{2}>0.
\]
Thus $p_*>0$ and $a<1/\rho$. As $s/n<\rho$,
\[
 p_s=1-\frac{as}{n}\ge1-a\rho=p_*\ge\frac{\alpha}{2}.
\]
It follows that
\[
 T_s=\frac{a}{p_s}\le T_*=
 \frac{a}{p_*}\le\frac{3/2}{\alpha/2}=\frac3\alpha.
\]
Therefore
\[
 T_*\xi\le\frac3\alpha\frac{\alpha^3}{100}
 =\frac{3\alpha^2}{100}
 <\frac{\alpha^2}{2}=(a-1)\alpha.
\]
Using $T_*=a/p_*$, this gives
\[
 \xi<\frac{a-1}{a}p_*\alpha<p_*\alpha<\alpha,
\]
since $(a-1)/a<1$ and $p_*<1$. Finally, $1/n<\xi/4$ gives
$\xi n>4$, so
\[
 q=\lfloor\xi n\rfloor\ge\xi n-1\ge\frac{\xi n}{2}\ge1.
\]
As $\xi<1$, also $q\le n$.
\end{proof}

For the rest of the construction, assume
\begin{equation}\label{eq:setup-size}
 n>T_*,\qquad 1/n<\xi/4.
\end{equation}

\begin{lemma}\label{lem:initial}
At stage $1$, for every $2\le j\le r$, every set $Q$ of $q$
labels, and every $W\subseteq V_j$,
\[
 \frac{1}{qn}\sum_{i\in Q}\deg_W(i)\ge \frac{|W|}{n}-\frac1n.
\]
In particular, \eqref{eq:check} holds at stage $1$.
\end{lemma}

\begin{proof}
At stage $1$, each label represents one vertex of $V_1$. The edges
between $V_1$ and $V_j$ form a matching, so each such vertex has at
most one neighbor in $V_j$. Fix $W\subseteq V_j$ and put
$\theta=|W|/n$. Each
label is compatible with at least $|W|-1$ vertices of $W$. By~\eqref{eq:setup-size} and~\eqref{eq:q-bounds},
we have $1\le q\le n$, so for each set $Q$ of
$q$ labels,
\[
 \frac{1}{qn}\sum_{i\in Q}\deg_W(i)
 \ge\frac{q(|W|-1)}{qn}
 =\theta-\frac1n.
\]
Since $1/n<\xi$, $p_1\le1$, and $\theta\ge0$,
\[
 \theta-\frac1n\ge\theta-\xi\ge p_1\theta-\xi.
\]
Taking $j=2$ gives \eqref{eq:check} for $s=1$.
\end{proof}

\section{Matchings in bipartite graphs}\label{sec:matching}

A \emph{fractional perfect matching} in a bipartite graph $B$ assigns a
nonnegative weight $f_{iv}$ to each edge $iv$ so that the weights on
the edges incident with each vertex sum to $1$. If each side of $B$
has $n$ vertices, then the total edge weight is $n$, since summing
the weights at the vertices on one side counts every edge exactly
once. An \emph{ordinary perfect matching} is the special case in which
each edge has weight $0$ or $1$.

\begin{lemma}\label{lem:cap}
If \eqref{eq:check} holds at stage $s$, then $B_s$ has a fractional perfect matching $f$ with $f_{iv}\le T_s/n$ on every edge.
\end{lemma}

\begin{proof}
Direct the edges of $B_s$ from left to right and give them capacity $T_s/n$. Add a source with an edge of capacity $1$ to each left vertex, and a sink with an edge of capacity $1$ from each right vertex. A \emph{flow} assigns nonnegative values to these edges, at most their capacities. At every vertex other than the source and sink, the incoming and outgoing totals must be equal. The value of the flow is the total leaving the source.

We seek a flow of value $n$. A \emph{cut} partitions the vertices into two sets, one containing the source and the other containing the sink. Its capacity is the sum of the capacities of the edges directed from the source side to the sink side. By the maximum-flow minimum-cut theorem, a flow of value $n$ exists if every cut has capacity at least $n$.

Consider any cut. Let $S$ be the set of left vertices on the source side, and let $W$ be the set of right vertices on the sink side. Thus the source side contains exactly the source, $S$, and the right vertices outside $W$. Write $m(S,W)$ for the number of edges of $B_s$ between $S$ and $W$.

There are three kinds of edges directed across this cut. The edges from the source to the left vertices outside $S$ contribute $n-|S|$. The edges from $S$ to $W$ contribute $(T_s/n)m(S,W)$. The edges from the right vertices outside $W$ to the sink contribute $n-|W|$. The cut therefore has capacity
\[
 (n-|S|)+\frac{T_s}{n}m(S,W)+(n-|W|).
\]
This is at least $n$ precisely when
\begin{equation}\label{eq:flow}
 \frac{T_s}{n}m(S,W)\ge h,\qquad h=|S|+|W|-n.
\end{equation}
Every cut is represented by such a pair $S,W$. If $h\le0$, then~\eqref{eq:flow} follows from $m(S,W)\ge0$, so assume $h>0$.

Put $d=n-s$. Every vertex of $B_s$ has degree at least $d$. 
Since $s<r\le\rho n$ and $\alpha=1-\rho$, the bound $\xi<\alpha$ in~\eqref{eq:slack-bounds}
gives
\[
 \frac dn=1-\frac sn\ge1-\rho=\alpha>\xi.
\]
Hence $d>\xi n\ge q$, since $q=\lfloor\xi n\rfloor$.
Since $a>1$, we have
\[
 p_s=1-\frac{as}{n}\le1-\frac{s}{n}=\frac dn.
\]
Using $T_s=a/p_s$ from~\eqref{eq:stage-parameters}, we obtain $T_sd/n\ge a>1$.

First suppose $|S|\le d$. Each label in $S$ has at least $d$ neighbors in total, and at most $n-|W|$ of them lie outside $W$. Its number of neighbors in $W$ is therefore at least
\[
 d-(n-|W|)=h+d-|S|.
\]
Summing over the labels in $S$ gives
\[
 \begin{aligned}
 m(S,W)&\ge |S|(h+d-|S|)\\
 &=dh+(d-|S|)(|S|-h)\\
 &\ge dh.
 \end{aligned}
\]
For the last step, $d-|S|\ge0$ by assumption, and $|S|-h=n-|W|\ge0$. It follows that
\[
 \frac{T_s}{n}m(S,W)\ge\frac{T_sd}{n}h\ge ah\ge h,
\]
as required. If $|W|\le d$, the same argument applies with the two sides exchanged, since the degree bound holds on both sides.

It remains to consider $|S|,|W|>d$. Put $\theta=|W|/n$. We have $|S|>q$, so we can apply~\eqref{eq:check} to every subset $Q\subseteq S$ of size $q$. There are $\binom{|S|}{q}$ subsets of $S$ of size $q$, and each
label belongs to exactly $\binom{|S|-1}{q-1}$ of them. Summing
\eqref{eq:check} over these subsets counts each term $\deg_W(i)$
exactly $\binom{|S|-1}{q-1}$ times. Since the right side is the
same for every $Q$, dividing by $\binom{|S|}{q}$ gives
\[
 \frac{\binom{|S|-1}{q-1}}{\binom{|S|}{q}}
 \frac{1}{qn}\sum_{i\in S}\deg_W(i)
 \ge p_s\theta-\xi.
\]
The ratio of binomial coefficients is $q/|S|$, and
$\sum_{i\in S}\deg_W(i)=m(S,W)$. Hence
\begin{equation}\label{eq:cut-edge-count}
 \frac{m(S,W)}{|S|n}\ge p_s\theta-\xi,
 \qquad\text{so}\qquad
 m(S,W)\ge |S|n(p_s\theta-\xi).
\end{equation}
Since $|W|>d$, we have $\theta>d/n\ge\alpha$. By~\eqref{eq:stage-bounds} and~\eqref{eq:slack-bounds},
\[
 T_s\xi\le T_*\xi<(a-1)\alpha<(a-1)\theta.
\]
Consequently $a\theta-T_s\xi>a\theta-(a-1)\theta=\theta$. Multiplying the lower bound for $m(S,W)$ in~\eqref{eq:cut-edge-count} by $T_s/n$ and using $T_sp_s=a$, we get
\[
 \begin{aligned}
 \frac{T_s}{n}m(S,W)
 &\ge |S|T_s(p_s\theta-\xi)\\
 &=|S|(a\theta-T_s\xi)\\
 &\ge |S|\theta=\frac{|S||W|}{n}.
 \end{aligned}
\]
Finally,
\[
 \frac{|S||W|}{n}-h
 =\frac{|S||W|-n|S|-n|W|+n^2}{n}
 =\frac{(n-|S|)(n-|W|)}{n}\ge0,
\]
because $|S|,|W|\le n$. This proves~\eqref{eq:flow} in the remaining case.

There is therefore a flow of value $n$. The $n$ edges from the source have capacity $1$ and total flow $n$, so each carries flow $1$. The same holds at the sink.
Set $f_{iv}$ equal to the flow on $iv$. Flow conservation gives
\[
 \sum_{v:\,iv\in E(B_s)}f_{iv}=1
 \quad\text{for each label }i,\qquad
 \sum_{i:\,iv\in E(B_s)}f_{iv}=1
 \quad\text{for each }v\in V_{s+1}.
\]
Since all edge flows are nonnegative, these weights form a
fractional perfect matching of $B_s$. The edge capacities give
$f_{iv}\le T_s/n$.
\end{proof}

An \emph{augmenting path} has unmatched endpoints and alternates between edges outside and inside the matching. Its first and last edges are outside the matching, so exchanging the two kinds of edges increases the matching size by one. The next proof finds an augmenting path of length at most five, which deletes at most two matching edges.

\begin{lemma}\label{lem:short}
Let $0\le\ell<d\le n$ be integers. Suppose a bipartite graph has two parts of size $n$, minimum degree at least $d$, and every set of $d$ right vertices has at least $n-\ell$ neighbors. Every matching $M$ of size $n-k$ admits a perfect matching $M^+$ with $|M\setminus M^+|\le2k$.
\end{lemma}

\begin{proof}
To enlarge a nonperfect matching $M$, choose an unmatched right vertex $b$. If $b$ has an unmatched left neighbor $c$, the edge $cb$ is an augmenting path of length one. Adding it enlarges $M$ without deleting any edge.

Otherwise every neighbor of $b$ is matched. Let $D$ be the set of right vertices matched to these neighbors. Different neighbors of $b$ are matched to different right vertices, so $|D|=\deg(b)\ge d$. Write $N(D)$ for the set of left vertices with a neighbor in $D$. Applying the hypothesis to any $d$ vertices of $D$ gives
\[
 |N(D)|\ge n-\ell.
\]
If $N(D)$ contains an unmatched vertex $c$, choose a neighbor $t\in D$ of $c$. Let $w$ be matched to $t$. By the definition of $D$, the vertex $w$ is a neighbor of $b$. Hence $c,t,w,b$ is an augmenting path of length three, with $tw$ its only matching edge.

It remains to consider the case in which every vertex of $N(D)$ is matched. Let $U$ be the set of right vertices matched to them. Then
\[
 |U|=|N(D)|\ge n-\ell.
\]
Also $D\subseteq U$. Indeed, for each $t\in D$, let $w$ be
the left vertex matched to $t$. Since $wt$ is an edge and
$t\in D$, we have $w\in N(D)$. Thus $t\in U$ by the
definition of $U$.

Choose an unmatched left vertex $c$, which exists because $M$
is not perfect and the two parts have the same size. At most
$\ell$ right vertices lie outside $U$, whereas $c$ has at least
$d>\ell$ neighbors. Thus $c$ has a right neighbor $z\in U$.
Since every vertex of $N(D)$ is matched whereas $c$ is unmatched,
we have $c\notin N(D)$. Hence $c$ has no neighbor in $D$.
As $cz$ is an edge, it follows that $z\in U\setminus D$.

Let $v$ be the left vertex matched to $z$. By the definition
of $U$, we have $v\in N(D)$, so $v$ has a neighbor $t\in D$.
By the definition of $D$, the vertex $t$ is matched to a left
neighbor $w$ of $b$. Consider the sequence
\[
 c,z,v,t,w,b.
\]
The edges $zv$ and $tw$ belong to $M$. The edges $cz$ and $wb$
do not belong to $M$, since $c$ and $b$ are unmatched.
Also $z\ne t$, because $z\notin D$ and $t\in D$. Since $v$
is matched to $z$, the edge $vt$ does not belong to $M$.
The edges therefore alternate between edges outside and inside $M$.

To verify that this sequence is a path, we also check that no
vertex is repeated. The right vertices $z$ and $t$ are distinct,
so their matched left vertices $v$ and $w$ are distinct as well.
The unmatched vertex $c$ differs from $v,w$, and the unmatched
vertex $b$ differs from $z,t$. Vertices on opposite sides are
distinct. Thus $c,z,v,t,w,b$ is an augmenting path of length five.

Starting with the given matching $M$ of size $n-k$, apply the
argument above whenever the current matching is not perfect.
Each augmentation removes at most two edges from the current
matching and inserts one more edge than it removes. Thus the
size of the matching increases by one. The graph does not change,
so its degree and neighborhood conditions remain valid.
The same argument can therefore be repeated until the matching
is perfect.

After $k$ augmentations, the resulting matching $M^+$ has size $n$
and is therefore perfect. Every edge of the original $M$ that is
absent from $M^+$ must have been removed during at least one
augmentation. Since each augmentation removes at most two edges,
\[
 |M\setminus M^+|\le 2k.
\]
If $k=0$, take $M^+=M$.
\end{proof}

\begin{lemma}\label{lem:completion}
If~\eqref{eq:check} holds at stage $s$, then $B_s$ satisfies
the hypotheses of Lemma~\ref{lem:short} with $d=n-s$ and $\ell=q$.
In particular, every matching $M$ of $B_s$ with $|M|=n-k$
admits a perfect matching $M^+$ such that
\[
 |M\setminus M^+|\le 2k.
\]
\end{lemma}

\begin{proof}
We verify the hypotheses of Lemma~\ref{lem:short}.
Take $d=n-s$ and $\ell=q=\lfloor\xi n\rfloor$. Every vertex has at most $s$ non-neighbors on the opposite side, so its degree is at least $n-s=d$. Since $s<r\le\rho n$ and $\xi<\alpha$ by~\eqref{eq:slack-bounds}, we have
\[
 \frac dn=1-\frac sn\ge1-\rho=\alpha>\xi.
\]
Thus $d>\xi n\ge q=\ell$, which verifies $0\le\ell<d\le n$.

Let $W$ be any set of $d$ right vertices. We now show that
$W$ has at least $n-q+1$ neighbors. Suppose, for a contradiction,
that at least $q$ labels have no neighbor in $W$, and choose
a set $Q$ of exactly $q$ such labels. Then $\deg_W(i)=0$
for every $i\in Q$. Applying~\eqref{eq:check}, with
$|W|/n=d/n$, gives
\[
 \begin{aligned}
 0=\frac{1}{qn}\sum_{i\in Q}\deg_W(i)
 &\ge p_s\frac dn-\xi\\
 &\ge p_*\alpha-\xi\\
 &>0.
 \end{aligned}
\]
The second inequality uses $d/n\ge\alpha$ and
$p_s\ge p_*>0$, the latter by~\eqref{eq:stage-bounds}.
The last inequality uses $\xi<p_*\alpha$, by
\eqref{eq:slack-bounds}.
This contradiction shows that at most $q-1$ labels have no
neighbor in $W$. Hence
\[
 |N(W)|\ge n-q+1\ge n-q=n-\ell.
\]
Together with the minimum degree bound and $0\le\ell<d\le n$
verified above, this establishes all the hypotheses of
Lemma~\ref{lem:short}. Applying that lemma gives the required
perfect matching $M^+$.
\end{proof}

\section{Random matchings}\label{sec:random}

The following result of Chekuri, Vondr\'ak and Zenklusen
\cite[Theorem 1.1]{CVZ} gives a random matching from a fractional
perfect matching. We state the result for bipartite graphs and use
it to bound the probability that the matching contains too few edges
from a fixed set.
For any edge set $F$, write $f(F)=\sum_{iv\in F}f_{iv}$.

\begin{lemma}[Chekuri, Vondr\'ak and Zenklusen]\label{lem:cvz}
Let $B$ be a bipartite graph with a fractional perfect matching $f$, and fix $0<\gamma\le1/2$. There is a random matching $M$ such that
\[
 \Prob(iv\in M)=(1-\gamma)f_{iv}
 \qquad\text{for every edge }iv.
\]
For every fixed edge set $F$, put $\mu_F=(1-\gamma)f(F)=\E|M\cap F|$. Then, for every $0<\zeta<1$,
\[
 \Prob\bigl(|M\cap F|<(1-\zeta)\mu_F\bigr)
 \le\exp(-\mu_F\gamma\zeta^2/20).
\]
The distribution of $M$ depends only on $B,f,\gamma$, and not on $F$ or $\zeta$.
\end{lemma}

We combine Lemmas~\ref{lem:short} and~\ref{lem:cvz} to obtain a perfect matching and bound the number of its edges in a fixed set.

\begin{lemma}\label{lem:round}
Suppose $B$ satisfies the hypotheses of Lemma~\ref{lem:short} and $f$ is a fractional perfect matching of $B$. For $0<\eta<1$, put $\kappa=(\eta/8)^3/40$. If $n\ge(\log2)/\kappa$, there is a random perfect matching $M^+$ of $B$ such that every fixed edge set $A$ satisfies
\[
 \Prob\bigl(|M^+\cap A|>f(A)+\eta n\bigr)\le\exp(-\kappa n).
\]
The distribution does not depend on $A$.
\end{lemma}

\begin{proof}
Set $\gamma=\eta/8$, so $0<\gamma<1/2$. Choose a random matching $M$ from Lemma~\ref{lem:cvz} and write $k=n-|M|$. 
Since $B$ satisfies the hypotheses of Lemma~\ref{lem:short},
for each possible matching $M$ we can choose a perfect matching
$M^+$ such that
\begin{equation}\label{eq:matching-changes}
 |M\setminus M^+|\le 2k=2(n-|M|).
\end{equation}
Fix this choice for each possible $M$ before considering any edge set $A$. This defines a distribution on perfect matchings that does not depend on $A$.

We first derive a probability estimate for any fixed edge set $F$. Put $X=|M\cap F|$ and $\mu=\E X=(1-\gamma)f(F)$. Summing the incident weights over one side of $B$ gives $f(E(B))=n$. Hence $f(F)\le n$ and $\mu\le n$. We claim that
\begin{equation}\label{eq:edge-count-tail}
 \Prob(X<\mu-\gamma n)\le\exp(-\gamma^3n/20).
\end{equation}
If $\mu\le\gamma n$, the event is empty because $X\ge0$. Otherwise set $\zeta=\gamma n/\mu$, which lies in $(0,1)$. Then $(1-\zeta)\mu=\mu-\gamma n$, and Lemma~\ref{lem:cvz} gives
\[
 \begin{aligned}
 \Prob(X<\mu-\gamma n)
 &\le\exp(-\mu\gamma\zeta^2/20)\\
 &=\exp\left(-\frac{\gamma^3n^2}{20\mu}\right)\\
 &\le\exp(-\gamma^3n/20).
 \end{aligned}
\]
The last inequality uses $\mu\le n$. This proves~\eqref{eq:edge-count-tail}.

Now fix $A\subseteq E(B)$ and put $C=E(B)\setminus A$. Applying~\eqref{eq:edge-count-tail} with $F=E(B)$, we have $X=|M|$ and $\mu=(1-\gamma)f(E(B))=(1-\gamma)n$. Except with probability at most $\exp(-\gamma^3n/20)$,
\begin{equation}\label{eq:matching-size}
 |M|\ge(1-\gamma)n-\gamma n=(1-2\gamma)n.
\end{equation}
Applying~\eqref{eq:edge-count-tail} with $F=C$, the mean is $(1-\gamma)f(C)$. Except with the same probability bound,
\begin{equation}\label{eq:complement-count}
 \begin{aligned}
 |M\cap C|&\ge(1-\gamma)f(C)-\gamma n\\
 &=f(C)-\gamma f(C)-\gamma n\\
 &\ge f(C)-2\gamma n,
 \end{aligned}
\end{equation}
where we used $f(C)\le n$. By adding the two failure probabilities, both estimates hold together except with probability at most $2\exp(-\gamma^3n/20)$.

When~\eqref{eq:matching-size} and~\eqref{eq:complement-count} hold, $k=n-|M|\le2\gamma n$. The bound on the number of deleted edges in~\eqref{eq:matching-changes} then gives
\[
 |M\setminus M^+|\le2k\le4\gamma n.
\]
At most this many edges of $M\cap C$ can be lost when we pass to $M^+$. Consequently,
\[
 \begin{aligned}
 |M^+\cap C|
 &\ge |M\cap C|-|M\setminus M^+|\\
 &\ge f(C)-2\gamma n-4\gamma n\\
 &=f(C)-6\gamma n.
 \end{aligned}
\]
Since $M^+$ has $n$ edges and $A,C$ partition $E(B)$, we obtain
\[
 \begin{aligned}
 |M^+\cap A|
 &=n-|M^+\cap C|\\
 &\le n-f(C)+6\gamma n\\
 &=f(A)+6\gamma n\\
 &<f(A)+\eta n.
 \end{aligned}
\]
Here $f(A)+f(C)=n$, and $6\gamma=3\eta/4<\eta$. 
This implies that the event
$|M^+\cap A|>f(A)+\eta n$ is contained in the union of the events
\[
 |M|<(1-2\gamma)n
 \qquad\text{and}\qquad
 |M\cap C|<f(C)-2\gamma n.
\]
Applying the union bound with the failure probabilities for~\eqref{eq:matching-size} and~\eqref{eq:complement-count}, we obtain
\[
 \begin{aligned}
 \Prob\bigl(|M^+\cap A|>f(A)+\eta n\bigr)
 &\le
 \Prob\bigl(|M|<(1-2\gamma)n\bigr)
 +\Prob\bigl(|M\cap C|<f(C)-2\gamma n\bigr)\\
 &\le 2\exp\left(-\frac{\gamma^3n}{20}\right)
 =\exp\left(\log 2-\frac{\gamma^3n}{20}\right)
 =\exp(\log 2-2\kappa n)\\
 &\le\exp(-\kappa n).
 \end{aligned}
\]
The last inequality follows from $n\ge(\log 2)/\kappa$.
\end{proof}

\section{Proof of Theorem~\ref{thm:main}}\label{sec:proof}

Set
\begin{equation}\label{eq:rounding-parameters}
 \eta=\frac{\xi^2}{8},\qquad
 \gamma=\frac{\eta}{8}=\frac{\xi^2}{64},\qquad
 \kappa=\frac{\gamma^3}{40},\qquad
 K=\left\lceil\frac{2\log4}{\kappa}\right\rceil.
\end{equation}
Assume also that
\begin{equation}\label{eq:size-conditions}
 n\ge\frac{\log2}{\kappa},\qquad n>\frac{4K}{3\xi}.
\end{equation}
The first inequality is the size assumption in Lemma~\ref{lem:round}.

We now specify the random assignments. Suppose that the process
has assigned the first $s$ parts. If~\eqref{eq:check} fails, stop.
Otherwise, choose a fractional perfect matching $f_s$ of $B_s$
with $(f_s)_{iv}\le T_s/n$ for every edge $iv$, as given by
Lemma~\ref{lem:cap}. By Lemma~\ref{lem:completion}, $B_s$
satisfies the hypotheses of Lemma~\ref{lem:short}. We can
therefore apply Lemma~\ref{lem:round} with this value of $\eta$
to choose a random perfect matching $M_s^+$. For each edge
$iv\in M_s^+$, add $v$ to the partial independent transversal
indexed by $i$. The choices of $f_s$ and of the distribution
of $M_s^+$ are fixed in advance for every possible assignment
of the first $s$ parts.

Since $\eta=\xi^2/8$ by~\eqref{eq:rounding-parameters} and $q\ge\xi n/2$ by~\eqref{eq:q-bounds}, the error term in Lemma~\ref{lem:round} satisfies
\begin{equation}\label{eq:eta-bound}
 \eta n=\frac{\xi^2n}{8}
 =\frac{\xi}{4}\frac{\xi n}{2}
 \le\frac{\xi q}{4}.
\end{equation}
The next bound applies to each choice of the parameters satisfying these conditions.

\begin{lemma}\label{lem:preserve}
Suppose $2\le r\le \rho n$, $n>T_*$, $1/n<\xi/4$, and~\eqref{eq:size-conditions} holds. The probability that the random process stops is at most
\begin{equation}\label{eq:stopping-bound}
 n^{K+1}\exp\bigl(-(\kappa K-\log4)n\bigr).
\end{equation}
If the process does not stop, it assigns all $r$ parts and produces a factor of independent transversals.
\end{lemma}

\begin{proof}
Fix a part $V_j$, where $2\le j\le r$, a set $Q$ of $q$ labels, and $W\subseteq V_j$. We track the compatible pairs in $Q\times W$ while the earlier parts $V_2,\ldots,V_{j-1}$ are assigned. Put $\theta=|W|/n$. At each reached stage $1\le s\le j-1$, define
\[
 \omega_s=\frac{1}{qn}\#\{(i,w)\in Q\times W:
          w\text{ is compatible with label }i\text{ at stage }s\}.
\]
Thus there are exactly $qn\omega_s$ compatible pairs, and $0\le \omega_s\le \theta$. Lemma~\ref{lem:initial} gives
\begin{equation}\label{eq:initial-density}
 \omega_1\ge\frac{q(|W|-1)}{qn}=\theta-\frac1n.
\end{equation}
Condition~\eqref{eq:check} just before assigning $V_j$ requires
$ \omega_{j-1}\ge p_{j-1}\theta-\xi$
for every set $Q$ of $q$ labels and every $W\subseteq V_j$.
For the fixed choice of $j,Q,W$, we bound the probability that
stage $j-1$ is reached and
$ \omega_{j-1}<p_{j-1}\theta-\xi.$

Consider an executed step $s\le j-2$. The matching $M_s^+$ assigns vertices $v\in V_{s+1}$, while $W\subseteq V_j$ lies in a later part. Let $A_s$ consist of the edges $iv\in E(B_s)$ with $i\in Q$ for which some $w\in W$ is compatible with $i$ at stage $s$ and satisfies $vw\in E(G)$. Thus $iv$ is an edge of the compatibility graph and $vw$ is an edge of the original graph.

Each edge $iv\in A_s$ corresponds to exactly one compatible pair $(i,w)$, since $v$ has at most one neighbor in $V_j$. Distinct edges correspond to distinct pairs, since $w$ has at most one neighbor in $V_{s+1}$. Using this correspondence and the edge weight bound $T_s/n$ from Lemma~\ref{lem:cap}, we obtain
\begin{equation}\label{eq:loss-weight}
 |A_s|\le qn\omega_s,
 \qquad
 f_s(A_s)\le\frac{T_s}{n}|A_s|\le qT_s\omega_s.
\end{equation}
Let $X_s=|M_s^+\cap A_s|$ count the edges of $A_s$ used in this assignment. If $(i,w)$ is currently compatible and $iv\in M_s^+$ assigns $v$ to label $i$, the pair becomes incompatible precisely when $vw\in E(G)$, because the earlier assignments remain unchanged. These lost pairs are exactly those corresponding to $M_s^+\cap A_s$, so their number is $X_s$. Adding a vertex cannot restore an incompatible pair. Consequently,
\[
 qn\omega_{s+1}=qn\omega_s-X_s.
\]
Dividing by $qn$ gives
\begin{equation}\label{eq:loss-count}
 \omega_{s+1}=\omega_s-\frac{X_s}{qn}.
\end{equation}

Condition on the assignments in the first $s$ parts. The graph $B_s$, the weights $f_s$, and the set $A_s$ are then fixed before $M_s^+$ is chosen. Lemma~\ref{lem:round}, applied with $A=A_s$, gives
\[
 \Prob\bigl(X_s>f_s(A_s)+\eta n
       \mid\text{assignments in the first }s\text{ parts}\bigr)
 \le\exp(-\kappa n).
\]
Adding the bounds $f_s(A_s)\le qT_s\omega_s$ from~\eqref{eq:loss-weight}
and $\eta n\le\xi q/4$ from~\eqref{eq:eta-bound} gives
\[
 f_s(A_s)+\eta n\le q(T_s\omega_s+\xi/4).
\]
If $X_s$ exceeds this larger threshold, it also exceeds
$f_s(A_s)+\eta n$. Thus the same probability bound gives
\begin{equation}\label{eq:conditional-loss}
 \Prob\bigl(X_s>q(T_s\omega_s+\xi/4)
       \mid\text{assignments in the first }s\text{ parts}\bigr)
 \le\exp(-\kappa n).
\end{equation}
The same $\kappa>0$ works for every possible set of previous assignments and every choice of $j,Q,W$, because it depends only on $\eta$.

For this fixed choice of $j,Q,W$, let $I_s=1$ if the event in~\eqref{eq:conditional-loss} occurs, and let $I_s=0$ otherwise. Set $I_s=0$ also if the step is not executed or $s\ge j-1$. The conditional probability bound still holds at unexecuted steps, where $I_s=0$. Their dependence on $j,Q,W$ is left out of the notation.

\begin{claim}\label{claim:bad-steps}
With probability at least
\[
 1-n^{K+1}\exp\bigl(-(\kappa K-\log4)n\bigr),
\]
every triple $(j,Q,W)$ has at most $K-1$ stages $s$ with $I_s=1$.
\end{claim}

\begin{proof}
For any prescribed stages $s_1<\cdots<s_K$, the earlier indicators are determined before the last random choice. Conditioning at the last stage and using~\eqref{eq:conditional-loss} gives
\[
 \begin{aligned}
 \Prob(I_{s_1}=\cdots=I_{s_K}=1)
 &\le e^{-\kappa n}\Prob(I_{s_1}=\cdots=I_{s_{K-1}}=1)\\
 &\le e^{-\kappa nK},
 \end{aligned}
\]
where the second line follows by repeating the same argument. No independence between the steps is needed.

The index $j$ ranges from $2$ to $r$, so it has $r-1$ possible
values. There are $\binom{n}{q}$ choices for $Q$, since it consists
of exactly $q$ of the $n$ labels. For each $j$, the set $W$ can be
any of the $2^n$ subsets of $V_j$. Since $r-1\le n$ and
$\binom{n}{q}\le 2^n$, the number of triples $(j,Q,W)$ is
\[
 (r-1)\binom{n}{q}2^n\le n\,2^n\,2^n=n4^n.
\]
Fix a triple $(j,Q,W)$. If $\sum_{s=1}^{j-2}I_s\ge K$, then
there are stages $1\le s_1<\cdots<s_K\le j-2$ such that
$I_{s_1}=\cdots=I_{s_K}=1$. Since $j-2<n$, there are at most
$n^K$ choices of these stages. For each choice, the probability
that all $K$ indicators equal $1$ is at most $e^{-\kappa nK}$.
The union bound therefore gives
\[
 \Prob\left(\sum_{s=1}^{j-2}I_s\ge K\right)
 \le n^K e^{-\kappa nK}.
\]
There are at most $n4^n$ choices of the triple $(j,Q,W)$.
Applying the union bound again, we obtain
\[
 \begin{aligned}
 &\Prob\left(\text{there is a triple }(j,Q,W)
       \text{ with }\sum_{s=1}^{j-2}I_s\ge K\right)\\
 &\hspace{15mm}\le n4^n\cdot n^K e^{-\kappa nK}\\
 &\hspace{15mm}=n^{K+1}\exp\bigl(-(\kappa K-\log4)n\bigr),
 \end{aligned}
\]
as desired.
\end{proof}

For each reached stage $s$, let
\[
 D_s=\max\{p_s\theta-\omega_s,0\}
\]
be the amount by which $\omega_s$ falls below $p_s\theta$.

\begin{claim}\label{claim:tracked-loss}
Fix a triple $(j,Q,W)$. If at most $K-1$ of its indicators $I_s$ equal $1$, then, whenever stage $j-1$ is reached,
\[
 D_{j-1}<\xi.
\]
\end{claim}

\begin{proof}
We estimate the total loss for this fixed triple. Consider an
executed step $s$. If $I_s=0$, then
$X_s\le q(T_s\omega_s+\xi/4)$ by the definition of $I_s$ in
terms of~\eqref{eq:conditional-loss}. Substituting this bound
in~\eqref{eq:loss-count} gives
\[
 \omega_{s+1}
 \ge\omega_s-\frac{T_s\omega_s+\xi/4}{n}
 =\left(1-\frac{T_s}{n}\right)\omega_s-\frac{\xi}{4n}.
\]

Now suppose that $I_s=1$. Since $M_s^+$ is a matching and every
edge of $A_s$ has its label in $Q$, we have $X_s\le q$.
By~\eqref{eq:loss-count},
\[
 \omega_{s+1}\ge\omega_s-\frac1n
 \ge\left(1-\frac{T_s}{n}\right)\omega_s
    -\frac{\xi/4+1}{n}.
\]
The second inequality follows from $T_s>0$, $\omega_s\ge0$,
and $\xi>0$. Both cases therefore give
\begin{equation}\label{eq:step}
 \omega_{s+1}\ge\left(1-\frac{T_s}{n}\right)\omega_s
 -\frac{\xi/4+I_s}{n}.
\end{equation}

By the definition of $p_s$ in~\eqref{eq:stage-parameters}, we have
$p_{s+1}=p_s-a/n$. The same equation gives $T_sp_s=a$. Hence
\[
 p_{s+1}
 =p_s-\frac{a}{n}
 =\left(1-\frac{T_s}{n}\right)p_s.
\]
Multiplying this identity by $\theta$ gives
\[
 p_{s+1}\theta
 =\left(1-\frac{T_s}{n}\right)p_s\theta.
\]
By~\eqref{eq:step}, we also have
\[
 -\omega_{s+1}
 \le-\left(1-\frac{T_s}{n}\right)\omega_s
    +\frac{\xi/4+I_s}{n}.
\]
Adding these two relations, we obtain
\[
 p_{s+1}\theta-\omega_{s+1}
 \le\left(1-\frac{T_s}{n}\right)(p_s\theta-\omega_s)
    +\frac{\xi/4+I_s}{n}.
\]

By~\eqref{eq:stage-bounds} and~\eqref{eq:setup-size}, we have
$0<T_s\le T_*<n$. Thus $0<1-T_s/n<1$. Since
$p_s\theta-\omega_s\le D_s$, it follows that
\[
 p_{s+1}\theta-\omega_{s+1}
 \le\left(1-\frac{T_s}{n}\right)D_s
       +\frac{\xi/4+I_s}{n}.
\]
The right side is nonnegative. By the definition of $D_{s+1}$,
\[
 \begin{aligned}
 D_{s+1}
 &=\max\{p_{s+1}\theta-\omega_{s+1},0\}\\
 &\le\max\left\{
 \left(1-\frac{T_s}{n}\right)D_s
 +\frac{\xi/4+I_s}{n},\,0\right\}\\
 &=\left(1-\frac{T_s}{n}\right)D_s
 +\frac{\xi/4+I_s}{n}.
 \end{aligned}
\]
Since $D_s\ge0$ and $1-T_s/n<1$, this gives
\begin{equation}\label{eq:deficit-step}
 D_{s+1}\le D_s+\frac{\xi/4+I_s}{n}.
\end{equation}

To bound $D_1$, we use~\eqref{eq:initial-density} and $p_1\le1$.
They give
\[
 p_1\theta-\omega_1
 \le(p_1-1)\theta+\frac1n
 \le\frac1n.
\]
Since $D_1=\max\{p_1\theta-\omega_1,0\}$, we have $D_1\le1/n$.

Suppose that stage $j-1$ is reached. Summing
\eqref{eq:deficit-step} over the preceding $j-2$ steps gives
\[
 D_{j-1}
 \le D_1+\frac{\xi(j-2)}{4n}
           +\frac1n\sum_{s=1}^{j-2}I_s.
\]
By assumption, at most $K-1$ of these indicators equal $1$.
Using $D_1\le1/n$ and $j-2<n$, we obtain
\[
 \begin{aligned}
 D_{j-1}
 &\le\frac1n+\frac{\xi(j-2)}{4n}+\frac{K-1}{n}\\
 &\le\frac{\xi}{4}+\frac{K}{n}\\
 &<\xi.
 \end{aligned}
\]
The last inequality follows from $n>4K/(3\xi)$ in
\eqref{eq:size-conditions}. When $j=2$, there are no preceding steps. In this case,
\[
 D_{j-1}=D_1\le\frac1n\le\frac Kn<\xi,
\]
where the last inequality follows from~\eqref{eq:size-conditions}.
\end{proof}

Suppose that the process stops just before assigning $V_j$.
The process stops only when condition~\eqref{eq:check} fails.
Thus there are a set $Q$ of $q$ labels and a set
$W\subseteq V_j$ such that
\[
 \omega_{j-1}<p_{j-1}\theta-\xi.
\]
By the definition of $D_{j-1}$, this gives $D_{j-1}>\xi$.

By Claim~\ref{claim:bad-steps}, with probability at least
\[
 1-n^{K+1}\exp\bigl(-(\kappa K-\log4)n\bigr),
\]
every triple $(j,Q,W)$ has at most $K-1$ indicators equal to
$1$. On this event, Claim~\ref{claim:tracked-loss} gives
$D_{j-1}<\xi$ for the triple above, a contradiction. Hence the
probability that the process stops is at most
\[
 n^{K+1}\exp\bigl(-(\kappa K-\log4)n\bigr).
\]

If the process does not stop, it assigns each part by a perfect
matching of its compatibility graph. The resulting independent
transversals therefore partition all vertices.

This completes the proof of Lemma~\ref{lem:preserve}.
\end{proof}

We now choose $\rho$ and verify all the conditions. Substituting $\xi=\alpha^3/100$ from~\eqref{eq:parameter-choice} into the formula for $\kappa$ in~\eqref{eq:rounding-parameters} gives
\begin{equation}\label{eq:kappa-alpha}
 \kappa=\frac{\xi^6}{40\cdot64^3}
  =\frac{\alpha^{18}}{40\cdot64^3\cdot100^6}.
\end{equation}
By~\eqref{eq:rounding-parameters}, $K\le 2(\log4)/\kappa+1$. Since $0<\alpha<1$, equation~\eqref{eq:kappa-alpha} gives an absolute constant $B>0$ such that
\begin{equation}\label{eq:K-bound}
 K\le B\alpha^{-18}.
\end{equation}
Choose an absolute constant $A\ge1$ with
\begin{equation}\label{eq:A-choice}
 A^{21}>400\cdot2^{21}B.
\end{equation}
For sufficiently large $r$, put
\begin{equation}\label{eq:rho-choice}
 \varepsilon=A r^{-1/21}\le1,\qquad
 \rho=\frac1{1+\varepsilon},\qquad
 \alpha=\frac{\varepsilon}{1+\varepsilon}.
\end{equation}
Then $\varepsilon/2\le\alpha\le\varepsilon$. For every $n\ge\lceil r+A r^{20/21}\rceil$, we have
\[
 n\ge(1+\varepsilon)r\ge r,\qquad r\le \rho n.
\]
Using $K\le B\alpha^{-18}$ from~\eqref{eq:K-bound} and $\xi=\alpha^3/100$ from~\eqref{eq:parameter-choice}, we obtain
\[
 \begin{aligned}
 \frac{K}{n\xi}
 &\le\frac{100B}{n\alpha^{21}}\\
 &\le\frac{100\cdot2^{21}B}{n\varepsilon^{21}}\\
 &=\frac{100\cdot2^{21}B}{A^{21}}\frac rn
 <\frac14.
 \end{aligned}
\]
The strict inequality uses the choice of $A$ in~\eqref{eq:A-choice} and $r/n\le1$. In particular, $n>4K/(3\xi)$. For this fixed $A$, we have $\alpha\ge A r^{-1/21}/2$ by~\eqref{eq:rho-choice}. Applying this lower bound to $T_*\le3/\alpha$ in~\eqref{eq:stage-bounds}, to $\xi=\alpha^3/100$, and to~\eqref{eq:kappa-alpha} gives
\[
 T_*\le\frac3{\alpha}=O(r^{1/21}),\qquad
 \frac1\xi=O(r^{1/7}),\qquad
 \frac1\kappa=O(r^{6/7}).
\]
Each is $o(r)$. Since $n\ge r$, all sufficiently large $r$ therefore satisfy $n>T_*$, $1/n<\xi/4$, and $n\ge(\log2)/\kappa$, for every allowed $n$. This verifies~\eqref{eq:setup-size} and~\eqref{eq:size-conditions}.

It remains to show that the upper bound~\eqref{eq:stopping-bound} on the probability of stopping is less than one. By~\eqref{eq:K-bound} and $\alpha\ge A r^{-1/21}/2$, we have $K=O(r^{6/7})$. The function $(\log \tau)/\tau$ decreases for $\tau\ge3$, so for $n\ge r\ge3$,
\[
 \frac{(K+1)\log n}{n}
 \le\frac{(K+1)\log r}{r}
 =O(r^{-1/7}\log r)\longrightarrow0.
\]
This estimate holds for all $n\ge r$. Since $K=\lceil 2(\log4)/\kappa\rceil$ by~\eqref{eq:rounding-parameters}, we have $\kappa K\ge2\log4$. Hence, for all sufficiently large $r$,
\[
 \begin{aligned}
 n^{K+1}\exp\bigl(-(\kappa K-\log4)n\bigr)
 &\le\exp\bigl((K+1)\log n-(\log4)n\bigr)\\
 &\le\exp\bigl(-\tfrac12(\log4)n\bigr)<1.
 \end{aligned}
\]
Thus the construction succeeds with positive probability for every $n\ge\lceil r+A r^{20/21}\rceil$. This proves the asserted existence of a factor. Combining this conclusion with the lower bound~\eqref{eq:lower-bound} gives
\[
 r\le f(r)\le\left\lceil r+A r^{20/21}\right\rceil
 \le r+A r^{20/21}+1.
\]
Hence $f(r)\le r+O(r^{20/21})$. Dividing by $r$ gives
\[
 1\le\frac{f(r)}r\le1+A r^{-1/21}+\frac1r\longrightarrow1,
\]
which completes the proof of Theorem~\ref{thm:main}.

\section*{Declaration of AI use}

The ideas underlying the proofs were developed independently by the author.
OpenAI tools were used solely to assist with language and presentation and to
discuss a few parts of the proofs. In the course of these discussions, the tools
suggested a separate, more involved argument for the stronger bound
$f(r)\le r+o((\log r)^2)$. The author instead chose to present the simpler
argument in this paper, which establishes the asymptotic form of the conjecture
with an explicit error bound. The author rewrote the text in his own words,
verified the mathematical arguments, and takes full responsibility for the final
manuscript.

\end{document}